\documentclass[onecolumn]{IEEEtran}

\usepackage{cite}
\usepackage{graphicx}
\usepackage{epsfig,psfrag}
\usepackage{balance}
\usepackage{array}
\usepackage{theorem}
\usepackage{amsmath}
\usepackage{amssymb}
\usepackage{flushend}

\begin{document}

\title{Discrete-Time Chaotic-Map Truly Random Number Generators: Design, Implementation, and Variability Analysis of the Zigzag Map}

\author{Hamid~Nejati, Ahmad~Beirami, and Warsame~H.~Ali

\thanks{H. Nejati is with the Department of Electrical Engineering and Computer Science, University of Michigan, Ann Arbor, MI, 48109 USA. e-mail:(hnejati@umich.edu)}

\thanks{A. Beirami is with the School of Electrical and Computer Engineering, Georgia Institute of Technology, Atlanta, GA, 30332 USA. e-mail:(beirami@ece.gatech.edu)}

\thanks{W. H. Ali is with the Electrical and Computer Engineering Department, Prairie View A\&M University, Prairie view, TX, 77446 USA. e-mail:(whali@pvamu.edu)}
}
\maketitle

\begin{abstract}
In this paper, we introduce a novel discrete chaotic map named zigzag map that demonstrates excellent chaotic behaviors and can be utilized in Truly Random Number Generators (TRNGs). We comprehensively investigate the map and explore its critical chaotic characteristics and parameters. We further present two circuit implementations for the zigzag map based on the switched current technique as well as the current-mode affine interpolation of the breakpoints. In practice, implementation variations can deteriorate the quality of the output sequence as a result of variation of the chaotic map parameters. In order to quantify the impact of variations on the map performance, we model the variations using a combination of theoretical analysis and  Monte-Carlo simulations on the circuits. We demonstrate that even in the presence of the map variations, a TRNG based on the zigzag map passes all of the NIST~$800-22$ statistical randomness tests using simple post processing of the output data.
\end{abstract}


\IEEEpeerreviewmaketitle

\section{Introduction}
\label{sec:intro}

Chaos is a non-periodic, long-term non-predictive behavior that can be generated by certain nonlinear dynamical systems~\cite{applied_chaos}. The chaotic systems are inherently deterministic given the initial state of the system. The chaotic behavior is a result of the exponential sensitivity of the system to the initial state that can not be exactly determined in practice. Chaotic waveforms have been extensively used in various research areas including biology for the modeling of the behavior of human organs~\cite{bio_chaos}, telecommunications for the modulation of the signals, chaotic oscillators, and chaotic encryption of the data~\cite{applied_chaos,MWSCASchaotic12,MWSCASchaotictune13}, and cryptography for the security of the system~\cite{crypto_handbook}.

Cryptography serves as the security block for the transmission of a message in the presence of eavesdroppers~\cite{crypto_handbook}. The security of a cryptographic system relies on the security of the Truly Random Number Generator (TRNG). A TRNG is capable of producing uncorrelated and unbiased binary sequences through a nondeterministic and irreproducible process~\cite{callegari05}.
\begin{figure} [t]
\vspace{0.05in}
\centering
\includegraphics[width=2.5in, angle=0]{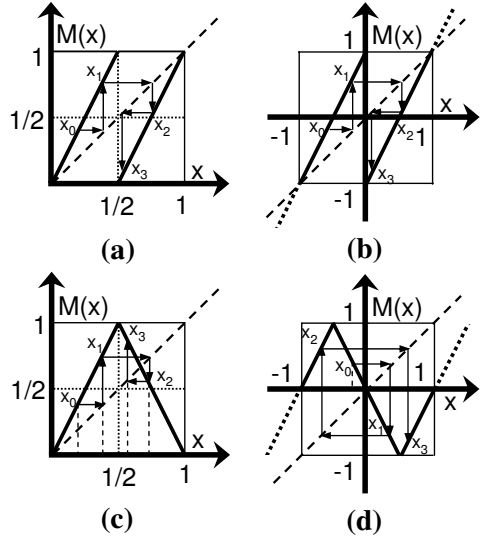}
\caption{(a)~The Bernoulli shift map. (b)~Practical implementation of the Bernoulli shift map using a pipelined ADC. (c)~The tent map. (d)~Our presented discrete chaotic map, called zigzag map.}
\label{fig:maps}
\end{figure}
Many practical cryptographic systems use a pseudo-random number generator (PRNG) as their random bit generator. PRNG is a deterministic system that produces a periodic limited-length binary sequences using a much shorter seed~\cite{crypto_handbook}. The period of the generated sequence grows exponentially with the length of the seed. The output of the PRNG with length much less than one period demonstrates random behavior in terms of bias and correlation between the generated bits. However, the entropy of any sequence generated using a PRNG is limited to the entropy of the seed sequence by the data processing theorem~\cite{Cover_book}. Consequently, a truly random seed maximizes the entropy of the sequence. Furthermore, the seed has to be periodically changed to keep the system secure. The large period of PRNG requires a large amount of storage capacity and redundant circuitry that makes PRNGs ineffective for high speed low cost cryptosystems~\cite{crypto_handbook}.

Software programs have been used for the implementation of TRNGs, such as time interval between keystrokes or mouse movements, system clock, content of I/O buffers, and operating system values such as network statistics~\cite{Hu07,crypto_handbook}.
Similar to most PRNGs, software based TRNGs generate output data using computer, which is in a large but finite number of states resulting in a pseudo-random output sequence~\cite{crypto_handbook}. In order to enhance the quality of these TRNGs, a complex combination of different sources of randomness should be used. This combination is practically inefficient, slow, and non-integrable.

On the other hand, hardware based TRNGs utilize inherent randomness of circuitry. The high entropy TRNG can be generated by the inherent randomness of the natural phenomena such as, the thermal and shot noise of the analog circuitry and the frequency instability (phase noise) of oscillators~\cite{crypto_handbook}.
Noise is presented in electronic circuits and it is also an important parameter in designing amplifiers and RF/analog systems~\cite{MWCASCS10,ISCASCS9,WOMICONmodel8,ALOGmodel7,ISCASMODELUWB6,GLSVLSICG4,ISCASCG3,JTCASLNA1}.
Direct amplification of the thermal noise of a resistor followed by a sampling stage seems to generate random sequences. However, this structure can not be used in integrated circuits due to the sensitivity to deterministic environmental variations such as power supply signals~\cite{petrie00}.
Most of the existing TRNGs are based on oscillator jitter, which is the variation of low-to-high and high-to-low transition points of the oscillator. If a low frequency (LF) oscillator with high jitter is used to sample the output of a high frequency (HF) oscillator, a random sequence can be obtained. Intel and Via use this structure in their security platforms~\cite{intel99,via03}.
However, the bit rate of a jitter-based random bit generator is limited to about $100$~kbps~\cite{liu05}, which is not suitable for high speed cryptosystems.
The use of continuous-time chaotic systems has been recently proposed for the generation of random numbers~\cite{yalcin04}. However, these methods suffer from the high correlation between the generated bits resulting in complicated post-processing of the generated bits. Therefore, these methods are not efficient in practice due to the sensitivity to the deterministic environmental variations that deteriorates the quality of the output random bit streams. In conclusion, these methods are not efficient in practice due
to the sensitivity to deterministic environmental variations that
deteriorates the quality of the output random bit streams.

A discrete-time chaotic map, formed by the iteration of the output value in a transformation function, can be used for the generation of random numbers. Simple piecewise affine input-output (I/O) characteristics have been extensively used for the generation of random bits, e.g., the Bernoulli map \cite{callegari05} and the tent map~\cite{bean94}. The entropy source of a chaotic map is the inherent noise of the system, which is amplified in the positive gain feedback loop by the iteration of the output signal in the map function~\cite{callegari05}. The output of the system will be unpredictable after a few first output samples. Pipelining multiple stages of the chaotic map circuit can increase the overall open loop gain of the system. The increased open loop gain can reduce the number of bits that need to be discarded in the beginning of the sequence due to the small inherent noise of the system. Pipelining multiple stages of the chaotic map circuit can also increase the speed of the whole system while contributing to the randomness and therefore improving the quality of the generated bits. The output of each stage is connected to the input of the successive stage. High Speed, capability of integration, and the high quality of the generated bits make the discrete-time chaotic maps excellent candidates for high speed embeddable random number generators.
Discrete-time chaotic TRNGs are much faster than the jitter based TRNGs. The high speed of these systems mainly owes to the circuit implementation and the pipelined structure, which can increase the output bit rate. The capability of being integrated is another important positive point for this kind of TRNGs.

Imperfections in practical implementation of TRNGs, so-called non-idealities, result in bias and correlation in the generated binary sequence of a truly random bit sequence. Hence, the post processing of the output data is required to improve the statistical characteristics of the bit sequence.  Von-Neumann~\cite{von-neumann51}, SHA-1~\cite{sha1}, and SHA-2~\cite{sha2} are the most widely used post-processing algorithms. All of these methods suffer from the trade off between bit-rate and the entropy of the generated binary sequence. Theory suggests that the knowledge of the non-ideality of the chaotic map can improve the efficiency of the post processing unit~\cite{ISIT11,ISIT12-gain,ISIT12-distributed}.
NIST~$800-22$ statistical test suite is a standard evaluation tool for the  examination of the quality of the generated bit stream~\cite{nist01}. This test can only prove that a large binary sequence is not random if it fails the tests. This test suite can not be used for the evaluation of the true randomness of the binary sequence. It is impossible to analytically prove the truly randomness of a sequence of data.

The practical application of both the tent map and the Bernoulli map can be hindered by noise and implementation errors, where they are unable to maintain the state of the system confined in the map~\cite{callegari05}. In the preliminary version of this work~\cite{Our_Midwest_Tent}
, we proposed zigzag map as a modification to the tent map that would make it robust to implementation variations and presented a very basic implementation.  The zigzag map is interchangeable with the tent map in practice while solving the confinement problem.  In this paper, we present two implementations for the zigzag map based on the switched current technique as well as the current-mode affine interpolation of the breakpoints for high speed applications. The current mode circuits are capable of producing high gain, less distortion, and high speed compared to voltage mode circuits. Less sensitivity to switching noise and process variations are other benefits of the current mode circuits. The effect of process variation on the parameters of the chaotic map is explored and a set of Monte Carlo simulations is used to quantize the map non-idealities. The second circuit demonstrated less sensitivity to process variations. The post-processing of the output data of a non-ideal map passes the randomness test suite.

This paper is organized as follows. In Section~\ref{sec:map}, the practical problems of the Bernoulli and tent maps are pointed out. Then, we present the zigzag map and investigate its chaotic characteristics. In Section~\ref{sec:circuit}, we explain the feasibility of implementation of the presented zigzag map. The deviation of the circuit I/O characteristic from the desired map deteriorates the quality of the output random bit sequences. This effect is modeled and analyzed in Section~\ref{sec:nonideality}. The generated sequence is fed to a post-processing unit and validated against NIST~$800-22$ randomness test suite. Section~\ref{sec:conclusion} concludes the paper.

\section{Discrete-Time Chaotic Map TRNGs}
\label{sec:map}

\subsection{Common Discrete-Time Chaotic Maps}
\label{subsec:discrete-time_maps}
Discrete-Time chaotic maps are a subclass of discrete-time nonlinear dynamical systems that can exhibit chaotic behaviors that are formed by the iteration of the output of the system into the system in every time step. A chaotic map is deterministic. If the initial state of a deterministic system is exactly known, the output behavior can be predicted exactly. Impossibility of exact knowledge of the initial state in practice is the reason of the chaotic behavior.

Bernoulli and tent maps are two most commonly used discrete-Time chaotic maps. The Bernoulli shift map, shown in Figure~\ref{fig:maps}(a), is a piecewise affine map that exhibits chaotic behavior. A $0$ is generated if $x_n<1/2$ and a $1$ is generated if $x_n>1/2$. The probability of all transitions are equal to $1/2$ regardless of the current state of the system. A pipelined ADC is used for practical implementation of the Bernoulli shift map. The I/O characteristic of a pipelined ADC is shown in Figure~\ref{fig:maps}(b). In practice, this map suffers from not being able to confine the output sequence in the chaotic region of $(-1,1)$. If for any reason the input value becomes slightly greater than $1$ or less than $-1$, the output stream gets out of the map and will never return back to the map. Decreasing the slope from $2$ and using  $1\frac{1}{2}$-bit ADC are the most widely solutions used to overcome this problem~\cite{li06,callegari05}, which results in very complicated structures.

The tent map, shown in Figure~\ref{fig:maps}(c), also exhibits chaotic behavior~\cite{bean94,morie00}. This system has been proved to be asymptotically stable with a uniform density distribution. $0$ is defined as $(0,1/2)$ and $1$ is defined as $(1/2,1)$. The proposed structures are sensitive to implementation variations~\cite{callegari97}. A tailed tent map was presented in \cite{callegari97}, which has a uniform density distribution. However, the process of generating bits from the tailed tent map is not straightforward. However, the post-processing required to compensate for non-idealities can result in complicated circuitry or reduced output bit-rate.
\begin{figure}[b]
\centering
\includegraphics[width=2.5in, angle=0]{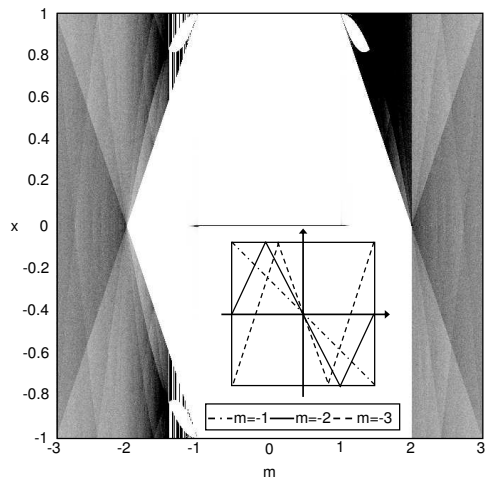}
\caption{Bifurcation diagram for the zigzag map with bifurcation parameter $m$, used in~(\ref{eq:bifurcation}). The generalized zigzag map with three different values of $m$ ($-1$, $-2$, and $-3$) is shown in the inset.}
\label{fig:bifurcation}
\end{figure}
\subsection{Presented Zigzag Map}
\label{subsec:presented_map}
We introduce  a novel discrete chaotic map and call it the zigzag map, shown in Figure~\ref{fig:maps}(d). It is notable that the output values of the zigzag map will alternate between positive and negative values. We will show later that the alternation between positive and negative values will be helpful in the presence of implementation non-idealities. As shown in Figure~\ref{fig:maps}(d), the tent map and the zigzag map have outputs that are equal in absolute value due to the symmetry but the output values alternate between positive and negative values in the zigzag map. Thus, the absolute value of the generated output stream is the same for both maps. In other words, if $x_1$, $x_2$ and $x_3$ are three successive outputs of the tent map, the outputs of the zigzag map will be given by $x_1$, $-x_2$ and $x_3$. Generation of bits from this map is straightforward. $|x_n|<1/2$ represents $0$ and $|x_n|>1/2$ represents $1$. The generalized zigzag map can be defined as
\begin{equation}
x_{n+1} = \left\{
\begin{array}{ll}
-m(x_n +\frac{2}{|m|}) & -1<x_n\leq-\frac{1}{|m|} \\
mx_n & -\frac{1}{|m|}<x_n\leq\frac{1}{|m|} \\
-m(x_n- \frac{2}{|m|}) & \frac{1}{|m|} <x_n \leq 1
\end{array} \right.,
\label{eq:bifurcation}
\end{equation}
where $m$ is a real number and $m \in (-3,3)$. In this equation, $m=-2$ represents the zigzag map. This formulation can be used to find the bifurcation diagram of the presented chaotic map with respect to the bifurcation parameter $m$ (Figure~\ref{fig:bifurcation}). A bifurcation diagram can be used to identify the chaotic behavior of a system. It can represent the possible steady state values of a dynamical system as a function of the bifurcation parameter.

Without loss of generality, the initial state of the system has been assumed to be a slightly positive number, i.e., $x_0 = 0^{+}$. We start analyzing the diagram from larger values of $m$. $m = 3$ is a critical point in the diagram. This case is at the verge of instability, since the system can be driven out of the map. In $2<m<3$, the system is chaotic and stable. Chaos can be observed for $1<m<2$. Since the initial state of the system has been assumed to be positive, the output value is confined to positive values. If the system gets out of the positive region and goes into the negative region it will get trapped in the negative region. For $|m|<1$, no chaotic behavior is observed. The output of the system will ultimately settle in zero. $-2<m<-1$ is also a chaotic region for the system. The output in this region alternates between positive and negative values. Note that there is no asymptotic density distribution since for odd and even time steps, the output values are in the positive and negative regions, respectively. In other words $\lim f_n(x)$ does not exist as $n \rightarrow \infty$. The bifurcation diagram of the tent map is the same as given here for $0<m<2$.

\section{Switched-Current Circuit Implementation for the Zigzag Map}
\label{sec:circuit}
\begin{figure}
\centering
\includegraphics[width=2.6in, angle=0]{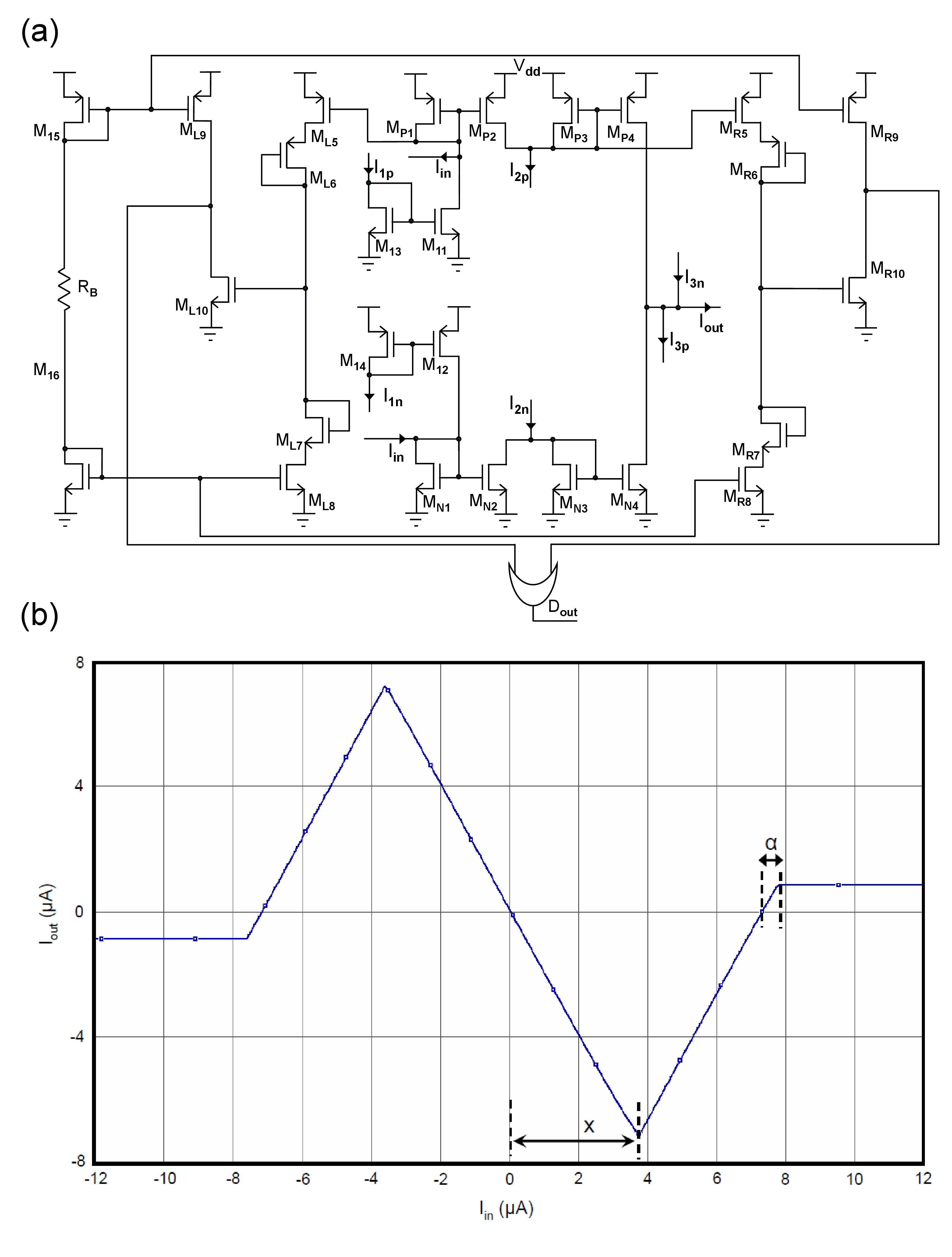}
\caption{(a)~The double breakpoint current mirrors for the generation of the I/O current characteristic of the zigzag map. (b)~The I/O current characteristic of the implemented circuit with the parameters ($x=3.6\mu A$ and $\alpha=0.5\mu A$).}
\label{fig:res_cir}
\end{figure}
Current mode circuits have recently attracted great attention~\cite{current_mode}. The signal is represented by current of branches instead of voltage of the nodes in the traditional voltage-mode circuits. Reduced distortion, high gain, low input and high output impedances, high speed, less sensitivity to switching noise and better ESD immunity are other excellent characteristics of current mode circuits that makes them one of the best candidates in the future telecommunications, analog signal processing, multiprocessors, and high speed computer interfaces~\cite{current_mode}.

In order to consider the effect of inter-die (different devices in one chip) or intra-die (among devices on various chips) variations, we design a current mode circuit with low sensitivity to the process variations. The symmetric double break point current mirrors are less sensitive to the process variations compared to other implementation of piecewise linear I/O current characteristics~\cite{NASA}. If we combine two different double break point circuits together, we can expect the I/O characteristic curve, shown in Figure~\ref{fig:maps}(d). One of the current mirrors is implemented by NMOS transistors, while for the sake of symmetry the other one is implemented by PMOS transistors. The double break point current mirrors do not accumulate errors originated from process variations. It is straightforward to design a desired piecewise linear I/O curve with this circuit. The preprocessing part of the circuit consists of two track and hold stages that enable the iteration of the output into the system.

The upper side of the circuit, shown in Figure~\ref{fig:res_cir}(a), can generate the double break point current mirror with the gain of $4$, while the gain of the lower part is $2$. Adding the output current of this two structures, we can generate the zigzag map by the I/O characteristic of the presented circuit. $I_{1}$, $I_{2}$, and $I_{3}$ are the current sources connected to the upper and lower side, discerning by an extra $p$ and $n$ index for upper (PMOS) and lower (NMOS) circuits, respectively. The values of these current sources are as follows:
\begin{eqnarray}
I_{1p}&=&x ,~I_{2p}=2x ,~I_{3p}=0,\nonumber\\
I_{1n}&=&2\alpha ,~I_{2n}=2x+\alpha ,~I_{3n}=2\alpha.
\end{eqnarray}

The I/O characteristics for $x=3.6\mu A$ and $\alpha=0.5\mu A$ as an example is shown in Figure~\ref{fig:res_cir}(b).

The reference current sources are implemented with resistors whose thermal noise can increase the randomness and uncertainty of the system.

Another common method to generate discrete chaotic maps is based on the current-mode affine interpolation of the breakpoints~\cite{758875}. This can be implemented by the affine interpolation of the value at the breakpoints~\cite{758875}. In this method, any affine partition of the map is implemented using one elementary block based on the detection of the breakpoints. Therefore, it is straight-forward to find out whether the output value is within a certain affine partition of the map.
\begin{figure}[b]
\centering
\includegraphics[width=2.0in, angle=-90]{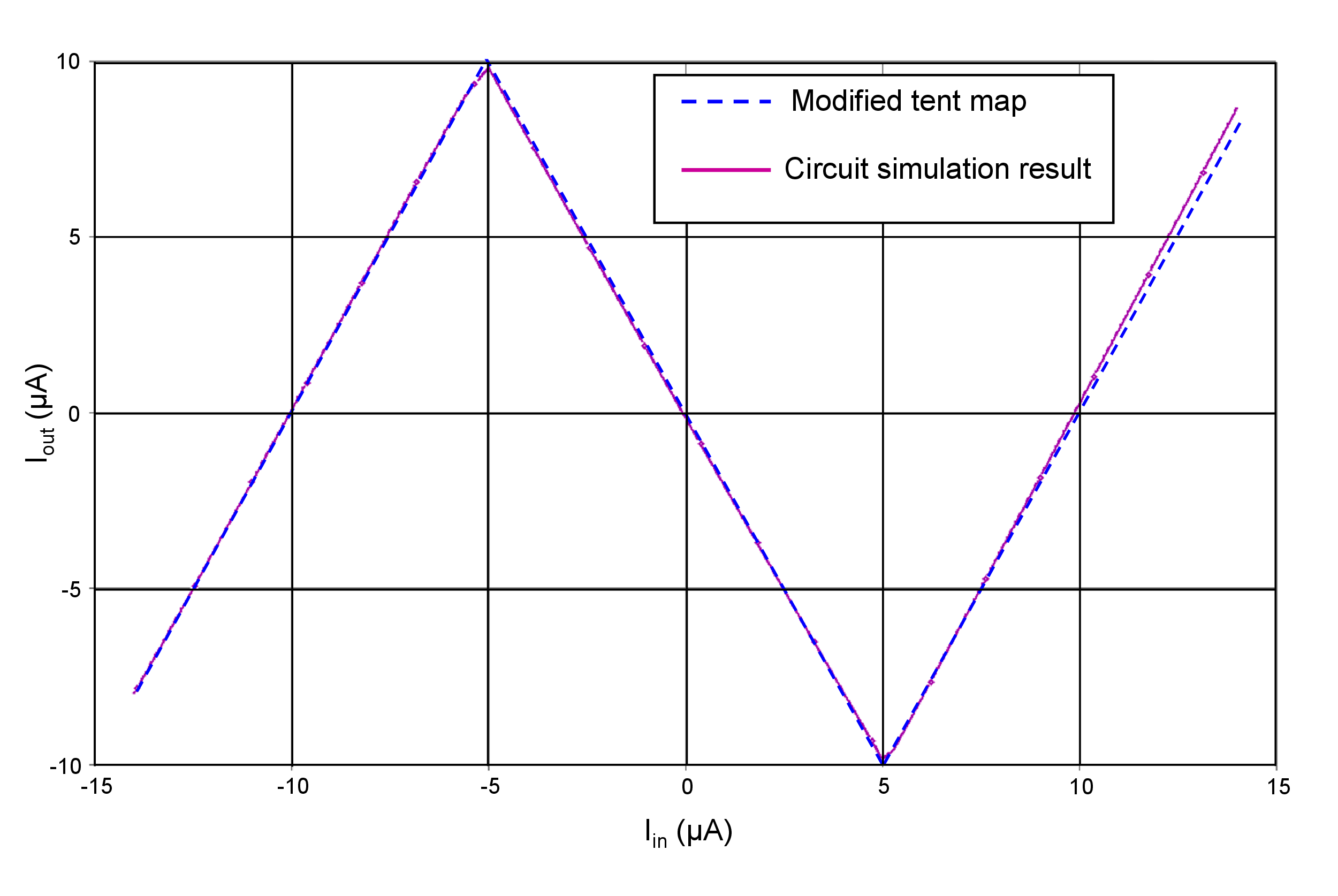}
\caption{The expected zigzag map (dashed line) and resulting I/O current characteristic of the zigzag map of the implemented circuit (solid line) based on the current-mode affine interpolation of the breakpoints with the parameters ($x=3.6\mu A$ and $\alpha=0.5\mu A$).}
\label{fig:result_circuit}
\end{figure}
In order to obtain the open-loop characteristic function of the circuit, the output of the circuit was connected to a matched load with the same impedance as the input impedance of the circuit. The simulations were carried out with CADENCE in IBM 0.13 $\mu m$ technology. In Figure~\ref{fig:res_cir}(b), the I/O characteristic of the first circuit obtained by simulation is shown.  $x$ and $\alpha$ are two parameters in the map, highlighted in Figure~\ref{fig:res_cir}(b). $\pm 2x$ is the required distance for the map. $\alpha$ is the extra length in the positive and negative side of the map. $\alpha=0$ can cause a trap point at zero, while considering non zero $\alpha$ can remove this state from the system and the possibility of trapping in the zero can be zero. The causes of slight changes in the simulation results will be taken into account in the next section. We performed HSPICE simulations in IBM 0.13 $\mu m$ technology for the second circuit. The resulting I/O current characteristics of the circuit is shown in Figure~\ref{fig:result_circuit}. The small variations of the map compared to the ideal map will be explained in the next section.

Four of the above maps are pipelined in a feedback loop. The noise floor of the analog circuits will set the initial state of the system that results in the non-predictive behavior of the system after a few first outputs. The inherent noise of the system is really small. However, first $m$ runs of the system amplify the noise to a detectable range. If the gain of the open loop system is $A$, we should discard first $m$ bit of the output stream given by $m=\log_A(P_d/N)$, where $P_d$ and $N$ are the detecting power and noise power, respectively. A single stage of the system has a gain $2$ and the whole system has a gain of $16$. The number of bits that need to be discarded from the beginning of the output sequence will be about $20$.

By investigating the transient response of both circuits with an input current pulse, we demonstrated that at the worst case the output would settle within $1\%$ of its final value in less than $12.5 ns$ and $15 ns$ for switched current circuit and interpolation of the values at breakpoints circuit, respectively.

\section{Modeling the Impact of Implementation Variations}
\label{sec:nonideality}
This section is divided into two parts. In sub section~\ref{subsec:nonideal_circuit}, we derive the impact of the variations of the physical parameters of the presented circuit as non-idealities of the chaotic map. In sub section~\ref{subsec:nonideal_map}, we investigate the transition matrix of the underlying Markov process in the presence of variations in the parameters of the chaotic map. We also quantify the bias and correlation in the output bit sequence.
\begin{figure}
\centering
\includegraphics[width=2.8in, angle=0]{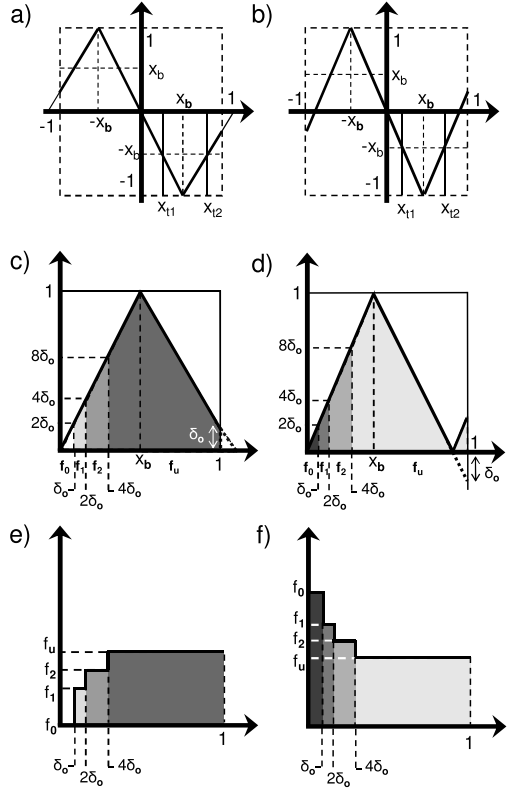}
\caption{The zigzag map with a fixed offset. The slope variations for
(a) negative and (b) positive changes. The equivalent tent map for the (a) and
(b) are shown in (c) and (d), respectively. The probability density function as
a function of the chaotic map parameter for (c) and (d) is depicted in (e) and
(f), respectively.}
\label{fig:nonideality_theory}
\end{figure}
\subsection{The Main Causes of Non-Ideality in the Presented Circuit}
\label{subsec:nonideal_circuit}
The circuit is based on the switched current technique. Hence, the process, voltage and temperature (PVT) variations mainly impact two different parts of the circuit, namely the sources of currents that are the current mirrors and the current limiter stages.

In a current mirror, the current $I$ that is supposed to mirror $I_s$, source current, is designed using the ratio of the width of transistors $W_2/W_1$. The current in a MOS transistor in saturation region is given by
\begin{equation}
I = \frac{1}{2} \mu C_{ox} \frac{W}{L} (V_{GS}-V_{th})^2(1 + \lambda V_{DS}).
\end{equation}
This suggests that
\begin{equation}
\frac{I_2}{I_1} = \frac{W_2}{W_1} \frac{L_1}{L_2} \frac{1 + \lambda_2 V_{DS2}}{1+\lambda_1 V_{DS1}}.
\end{equation}
Since we are not designing in deep sub micron processes, we can neglect the effect of variations in $\lambda$. We can approximate the variation in different process and $V_{DS}$ provided that the changes compared to the parameters are negligible. Hence,
\begin{equation}
\frac{I_2}{I_1}\hspace{-0.01in} \approx\hspace{-0.01in} \frac{W_2}{W_1} \hspace{-0.01in}(\hspace{-0.01in}1+2|\Delta W|\hspace{-0.01in}+2|\Delta L|\hspace{-0.01in}
+4\frac{|\Delta V_{th}|}{V_{gs}-V_{th}}+2\frac{\lambda |\Delta V_{DS}|}{1+\lambda V_{DS}})\hspace{-0.02in},
\end{equation}
where $\Delta W$, $\Delta L$, $\Delta V_{th}$, and $\Delta V_{DS}$ are the variation in width, length, threshold voltage, and the drain-source voltage of the transistors, respectively. If the $V_{DS}$ changes are large, we should compensate that by further tuning of the width of the transistors. We need to investigate the worst-case variations.
The impact of variations in $V_{DS}$ and $V_{th}$ is in the range of $1\%$ which can be ignored compared to other parameters~\cite{fiez91,ITRS}. Thus, the inner die variations for width of two different transistors in the current mirror will play the dominant role in altering the gain of the current mirror topology. The standard deviation of the variations are around $5\%$~\cite{nassif01,ITRS,ISQEDYIELD5,JALOGVAR2}.

The current limiter can not exactly provide the desired fixed current due to  the output resistance of the transistors. This non-ideality results in a small variation in the slope in the I/O characteristic of the circuit, which is proportional to $\lambda$ and the bias current of the circuit. While we can reduce this effect by further tuning of the width of transistors, the variation is inevitable.

\subsection{Investigating the Effect of Map Non-Idealities on Bias and Correlation of Output Data Sequence}
\label{subsec:nonideal_map}
\begin{figure}
\centering
\includegraphics[height=2.8in, angle=-90]{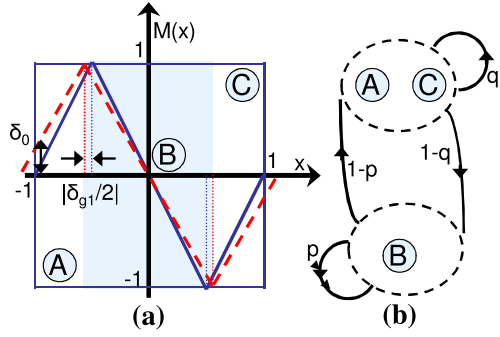}
\caption{a) The impact of variations on the zigzag map characteristics. b) The states of the Markov chain for truly random number generation.}
\label{fig:nonideal_map}
\end{figure}

\begin{figure*}
\vspace{0.05in}
\centering
\includegraphics[width=6.0in, angle=0]{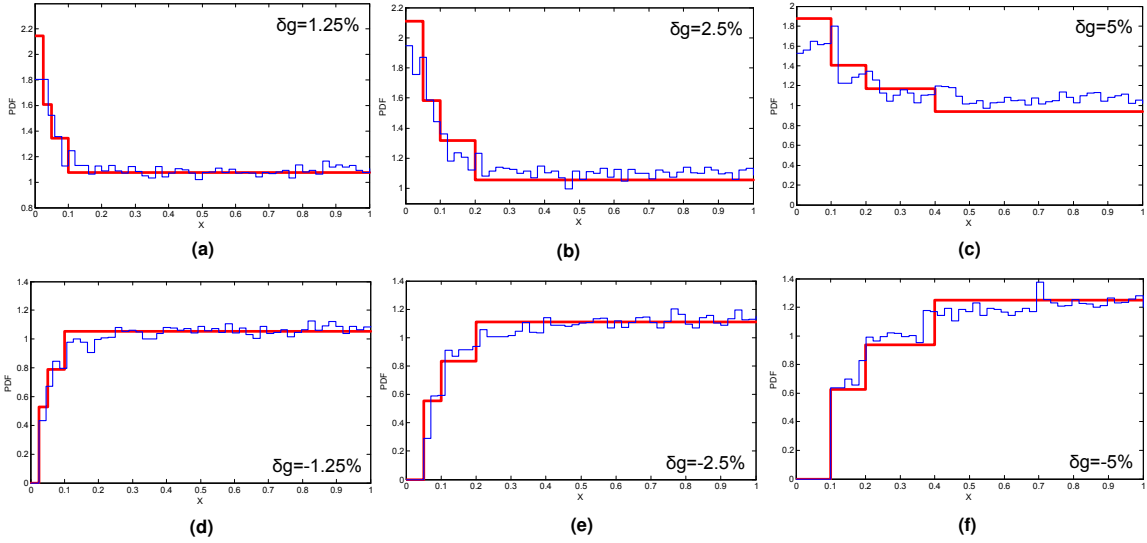}
\caption{The probability distribution function of the non-ideal system for different slope variations of (a) $+1.25\%$, (b)  $+2.5\%$, (c)  $+5\%$, (d)  $-1.25\%$, (e)  $-2.5\%$, and (f)  $-5\%$. The thick solid line demonstrates the presented model and the thin solid line represents the simulation results. }
\label{fig:nonideality_practice}
\end{figure*}

We first consider a map with symmetric slope variation as in Figures~\ref{fig:nonideality_theory}(a,b). We later extend to the non-symmetric case when we study the impact of post-processing.
 As discussed earlier, the ideal zigzag map has the same functionality as the tent map except for the alternations of the output between positive and negative values. Due to the symmetry in the map function and the bit generation process, the non-ideal maps can be analyzed using the maps shown in Figure~\ref{fig:nonideality_theory}(c,d). Since the latter are easier to work with, we use them for the rest of the derivations. Note that without loss of generality, every chaotic map can be scaled to $(-1,1)$. We assume the slope (gain) in the rising part of the map to be $2(1+\delta_{g1})$ and the gain in the falling part to be $-2(1+\delta_{g2})$.
The non-ideal map function $N(x): (0,1)\rightarrow (0,1)$ is thus given by
\begin{equation}
N(x) = \left\{
\begin{array}{ll}
2(1+\delta_{g1})x & x< x_b\\
1-2(1+\delta_{g2})(x-x_b) & x>x_b
\end{array}
\right. .
\end{equation}
where $x_b$ is abscissa of the boundary between rising and falling sections and is
\begin{equation}
x_b = \frac{1}{2(1+\delta_{g1})}.
\end{equation}
Note that our analysis of the symmetric non-ideal map is also valid for the non-ideal tent map. In this map, bit $0$ or $1$ is generated if $x\in (0,x_b)$ or $x \in (x_b,1)$, respectively. Throughout our analysis, we assume that variations are small, and we use first-order approximations for $\delta_{g1}$ and $\delta_{g2}$. This assumption is valid according to our analysis of the process variations in the current mirrors. Accordingly, $x_b$ can be easily shown to be
\begin{equation}
x_b \approx \frac{1}{2} (1-\delta_{g1}).
\end{equation}
Let $\delta_o$ be the ordinate of the end point of the map with abscissa $x=1$, which is given by
\begin{equation}
\delta_{o} = - ( \delta{g1} + \delta{g2}).
\end{equation}
The impact of $x_b$ and $\delta_0$ is illustrated in Figure~\ref{fig:nonideal_map}(a). The analysis of the two cases of $\delta_o >0$ and $\delta_o <0$ is slightly different. First, consider Figure~\ref{fig:nonideality_theory}(c), i.e., $\delta_o >0$. We want to find a good approximation for the asymptotic probability distribution $f(x)$ in the map. We approximate $f(x)$ using a four step distribution, as shown in Figure~\ref{fig:nonideality_theory}(e). Let $f_u$ denotes the approximate probability density at $x=1$.
According to the Frobenius-Perron operator, the distribution in $(0,\delta_{o})$ can be given by
\begin{equation}
f_0 = \frac{1}{2} f_0 = 0.
\end{equation}
We can calculate the distribution in $(\delta_o, 2\delta_o)$, $f_1$, as
\begin{equation}
f_1 = \frac{1}{2} f_0  + \frac{1}{2} f_u = \frac{1}{2} f_u.
\end{equation}
Finally, the probability distribution for $(2\delta_o, 4\delta_o)$, $f_2$, can also be calculated using the same approach to be
\begin{equation}
f_2 = \frac{1}{2} f_1 + \frac{1}{2} f_u = \frac{3}{4} f_u.
\end{equation}
The probability distribution in $(4\delta_o,1)$ is approximated using $f_u$. Now, $f_u$ can be calculated by noting that the probability distribution integrates to unity, i.e., $\int{f(x)dx}=1$. Thus,
\begin{equation}
f_u  = \frac{1}{1-2\delta_o} \approx 1 + 2\delta_o.
\end{equation}
The results of simulation are provided in Figure~\ref{fig:nonideality_practice}(a-c) for different values of $\delta_o$.
We can see that our approximate model captures the impact of the variations even for fairly large values of $\delta_o$.

The case of Figure~\ref{fig:nonideality_theory}(d), i.e., $\delta_o < 0$, can be analyzed similarly. The four step density distribution is shown is Figure~\ref{fig:nonideality_theory}(f). We derive the distribution in $(0,|\delta_o|)$ as
\begin{equation}
f_0 = \frac{1}{2} f_0 + \frac{1}{2} f_u + \frac{1}{2} f_u = 2f_u.
\end{equation}
The probability distribution in $(|\delta_o|,2|\delta_o|)$, i.e., $f_0$, the probability distribution in $(2|\delta_o|,4|\delta_o|)$, i.e., $f_1$, and  the distribution in $(|4\delta_o|,1)$, i.e., $f_2$, can be calculated in a similar way. W be $f_1 = \frac{3}{2} f_u$ and $f_2=\frac{5}{4}f_u$.
Finally, we have
\begin{equation}
f_u = \frac{1}{2(1-\delta_o)} \approx 1+2\delta_o.
\end{equation}
In the case of $\delta_o <0$, we have $f_u<1$, i.e., the value of $f(x)$ for $x \approx 1$ is less than unity. Again, the simulation results presented in Figure~\ref{fig:nonideality_practice}(d-f) show good agreement between the presented model and simulation results.


We model the bit generation process using a Markov chain as demonstrated in Figure~\ref{fig:nonideal_map}(b).
The transition probabilities of the map can be derived by integration of the probability distribution of the map to form the transition matrix of the Markov chain as
\begin{eqnarray}
P(0|0) = \frac{\int_{0}^{x_{t1}} {f(x)dx}}{\int_{0}^{x_b} {f(x)dx}}, & P(1|0) = \frac{\int_{x_{t1}}^{x_b} {f(x)dx}}{\int_{0}^{x_b} {f(x)dx}},\nonumber\\
P(0|1) = \frac{\int_{x_{t2}}^{1} {f(x)dx}}{\int_{x_b}^{1} {f(x)dx}}, &
P(1|1) = \frac{\int_{x_b}^{x_{t2}} {f(x)dx}}{\int_{x_b}^{1} {f(x)dx}},
\label{eq:transition_matrix}
\end{eqnarray}
where $x_{t1}$ and $x_{t2}$ are the points whose ordinates are equal to $x_b$ and are shown in Figure~\ref{fig:nonideality_theory}(a,b).

We use our approximate model for the evaluation of the integrals. First, we derive $p$ and $q$. i.e., $P(0|0)$ and $P(1|1)$. For fairly small variations ($|\delta_o| < 1/16$), the three regions $f_0$, $f_1$ and $f_2$ lie in $(0,x_{t1})$. This is indeed the region of interest for our purposes. We derive the transition probabilities as functions of $\delta_{g1}$ and $\delta_{g2}$. $p$ and $q$ are given by
\begin{eqnarray}
&&p = P(0|0) = \frac{1}{2} + \frac{3}{2} \delta_{g1} + 2 \delta_{g2}, \nonumber\\
&&q = P(1|1) = \frac{1}{2} - \frac{1}{2} \delta_{g2}.
\label{eq:p_q_1_1}
\end{eqnarray}

Thus far, we derived the transition probabilities of the underlying Markov chain in~(\ref{eq:transition_matrix}). Next, we demonstrate that these can actually be used to derive the bias and correlation in the
output bit sequence of the truly random number generator~\cite{MWSCASTRNG211}. It is straightforward to show that the bias is
\begin{equation}
b = \left|\frac{1}{2} - P(0)\right| = \left|\frac{1}{2} - P(1)\right| = \frac{|p-q|}{2 - p - q}.
\end{equation}
Further, $\lambda_1 = |p+q -1 |$ is the second eigenvalue of the transition Matrix. It can be shown that the correlation between the bits in the output sequence decays exponentially with their distance as $2^{-cn}$, where $c = -\log_2(\lambda_1)$.
Note that if $\delta_{g1} = \delta_{g2} = 0$, we will have $p = q = \frac{1}{2}$ and thus there is no bias ($b=0$) and no correlation ($c=\infty$) in the output sequence.


\section{Post processing}
\label{sec:postprocess}
We performed an extensive set of Monte Carlo simulations with HSPICE to determine the impact of process variations on the slope of I/O characteristic curve. The absolute value of the slope, characterized by $m$, should be 2 for the presented chaotic map. Monte Carlo simulation technique uses different randomly chosen values for the transistor parameters with a defined Gaussian distribution. We assume that the transistor parameter values are normally distributed and uncorrelated with a standard deviation of 2 percent, which is within the range of predicted values for current and future process technologies. Monte Carlo simulations demonstrate that the second circuit is less sensitive to process variations; therefore we use this circuit for the rest of this paper. By inspecting the variation of different linear pieces of the I/O characteristic curve as the objective of this Monte Carlo simulation, we investigate that the variation in slopes are less than 4 percent. To be more accurate, the standard deviation of the slopes are less than 2, 4, and 1 percent for the first, second, and third linear pieces, respectively.

Due to imperfections in practice, the existence of correlation and bias is inevitable that can be removed with post-processing unit. Von-Neumann, SHA-1, and SHA-2 post-processing algorithms are the most widely algorithms in practice and have proved to be effective in removal of the bias and correlation~\cite{zhou06a,wang05b,yalcin04}. In Von-Neumann algorithm, the sequence is divided into the blocks of length two. Every $01$ and $10$ are considered $0$ and $1$, respectively, while $00$ and $11$ are discarded. This results in a variable bit-rate that is not desired. SHA-1 and SHA-2 algorithms are proper for removing the correlation and bias from the bit streams. Although these methods are effective, the implementation requires FPGA or complex hardware configurations. The sequences of data should pass through the post-processing stage mainly because of inherent non-idealitiy factors in process or circuit.

\begin{figure}
\vspace{0.05in}
\centering
\includegraphics[width=2.8in, angle=0]{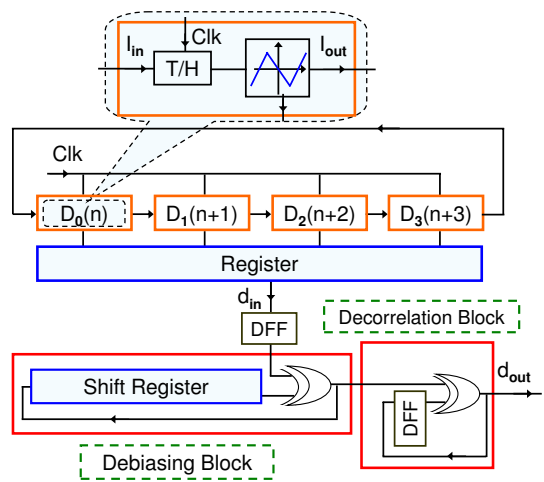}
\caption{The overall implementation of the presented TRNG with four pipelined chaotic maps. The post processing unit can remove the correlation and bias of the bit shuffled streams with decorrelation and debiasing blocks, respectively. }
\label{fig:post_processing}
\end{figure}

Our post-processing module is based on the method proposed by Addabbo et. al.~\cite{addabbo06b}. The method relies on two main facts. First, bits that are located further apart from each other have less correlation. The correlation approaches $0$ as the distance of the bits approaches infinity. The other fact is that XOR of two \emph{independent} biased sequence will result in a sequence with less bias. See~\cite{sunar07} for details. The debiasing block, shown in Figure~\ref{fig:post_processing}, consists of a shift register of length $l$ and an XOR gate. $z_{n-l}$, which is the output of time $n-l$ is XORed with $d_{in}(n)$ which is the input of time $n$. $l$ must be large enough so that $z_n$ and $z_{n+l}$ are almost uncorrelated. In practice, $l$ can be efficiently chosen by using the estimations made in the last sections. This has been proved to remove the bias of the output sequence in steady state. The decorrelation block is used to remove the correlation of the bits with a distance of $kl$. Note that we take the number of stages of shift register $l$ and the number of stages of the map $m$ to be prime with respect to each other, $(m,l) =1$. This will reduce greatly the correlation between bits of distance $kl$, where $k$ is not a multiplicand of $l$. The presented TRNG is robust to process and circuit variations, which provides simple post-processing stage compared to other chaotic maps used in TRNGs.

The output bit streams of the ideal and non-ideal I/O characteristics are tested
in NIST~$800-22$ test suite and the results are presented in Table~I. As listed in the table, before post-processing, the generated bits fail some of the tests but after passing through post-processing unit the generated bits pass all the tests.

In the generation of the random bits, we consider the worst case scenario, which consists of four stages with similar $I_{out}/I_{in}$ characteristics. In practice, these four stages are not exactly the same, therefore, the mismatch can enhance the quality of random output sequence. Utilization of four stages can also increase the total bit rate by multiplying the internal clock (track and hold clock) of the system by the number of cascading stages. The output bits of all stages go to a 4 bit register, known as bit shuffling~\cite{Gerosa02}. The uniformity of the output bit stream is more than the inputs. The output of the register feed to the post processing unit.

The presented chaotic map can generate random bits, which have lower sensitivity
to the circuit parameter fluctuations caused by process variation or circuit element mismatches. In addition to low power consumption ($1.42~mW$), the presented circuit can operate in clock frequencies as high as $10~MHz$ with the overall bit rate of $40~Mbps$. To calculate the worst case power consumption of the circuit, we calculate the average power consumption of the circuit with different input currents. Due to randomness and chaotic behavior of the presented chaotic map, the probability of the system to be in any point in the steady state is the same.

\begin{table}
\begin{center}
\vspace{0.3in}
\caption{NIST-$800-22$ test suite results for the non-ideal zigzag map before and after applying the post-processing unit}
\small
\vspace{0.1in}
\renewcommand{\arraystretch}{1.0}
\small
\textwidth=0.5in
\begin{tabular}{|c|c|c|c|}
\hline
 NIST$-800-22$ & Nonideal map w & Nonideal map w/o\\
  test suite  & post-processing & post-processing\\
\hline
apen & 0.513673 & 0.431503 \\
\hline
 Block frequency & 0.410214 & 0.94278 \\
\hline
 Cumulative sums & 0.80785 & 0.095407 \\
\hline
 Cumulative sums & 0.975089 & 0.03606 \\
\hline
FFT & 0.871131 & 0.721176  \\
\hline
Frequency & 0.575479 & 0.047704  \\
\hline
Linear complexity & 0.082259 & 0.017841  \\
\hline
 Longest run & 0.024845 & 0.553323  \\
\hline
 Non-periodic templates & 100$\%$ Success & 99.32$\%$ Success  \\
\hline
 Overlapping templates & 0.736301 & 0.510453  \\
\hline
Rank & 0.98771 & 0.445713  \\
\hline
Runs & 0.044596 & 0.689421  \\
\hline
 Serial & 0.443073 & Failure  \\
\hline
 Serial & 0.548205 & Failure  \\
\hline
bias & 0.04$\%$ & 0.50$\%$  \\
\hline
\end{tabular}
\end{center}
\label{tab:NIST2}
\end{table}
\section{Conclusion}
\label{sec:conclusion}
In this paper, we presented a TRNG based on a novel discrete-time chaotic map, namely the zigzag map, that is capable of being used in high speed embedded cryptographic systems. The proposed map proved to be chaotic and has been implemented by switched current technology. The random bits are generated due to the transition among different points of the chaotic map. Pipelining four stages in a positive feedback loop, we can amplify the inherent noise of the system within a detectable range, which results in random number generation. A modeling methodology has been developed to predict the performance of the system in the presence of non-idealities. Non-ideality sources are modeled and their effect on the performance of the output binary sequence (correlation and bias) has been studied. The non-ideal binary stream has been fed to a post processing module to increase the quality of the random digits. The presented system has been demonstrated to generate random bits that pass the NIST~$800-22$ test suite.

\section*{Acknowledgements}

This work was done in part when the authors were affiliated with Rice University. The authors would like to thank Professor Yehia Massoud for the valuable discussions, comments and suggestions that helped to improve the quality of this paper.

\end{document}